\documentclass[12pt]{article}
\usepackage{amsmath,amsfonts,amssymb,amsthm,enumerate,longtable}

\newcommand{\R}{{\mathbb R}}
\newcommand{\C}{{\mathbb C}}
\newcommand{\Z}{{\mathbb Z}}

\newcommand{\1}{{\bf 1}}
\newcommand{\F}{{\mathbb F}}

\newcommand{\BW}{{\rm BW}_{32}}
\newcommand{\sBW}{{\rm BW}_{16}}

\newcommand{\Aut}{{\rm Aut\, }}

\newcommand{\w}{\omega}
\newcommand{\RM}{{\rm RM}}
\newcommand{\eqa}{\begin{eqnarray}}
\newcommand{\eeqa}{\end{eqnarray}}
\newcommand{\eqn}{\begin{eqnarray*}}
\newcommand{\eeqn}{\end{eqnarray*}}

\newcommand{\Hom}{{\rm Hom}_\Z}
\newcommand{\Irr}{{\rm Irr}}

\newtheorem{dfn}{Definition}[section]
\newtheorem{pro}[dfn]{Proposition}
\newtheorem{thm}[dfn]{Theorem}

\newtheorem{lem}[dfn]{Lemma}
\newtheorem{cor}[dfn]{Corollary}
\newtheorem{rem}[dfn]{Remark}
\newtheorem{note}[dfn]{Note}
\newtheorem{example}[dfn]{Example}

\newcommand{\Ker}{{\rm Ker\ }}
\newcommand{\NO}{\,{\raise0.25em\hbox{$\mathop{\hphantom{\cdot}}\limits^{_{\circ}}_{^{\circ}}$}}\,}
\newcommand{\IntTK}[3]{
\scriptstyle\hphantom{\scriptstyle#2}\hfil#1\hfil\hphantom{\scriptstyle#3}%
\atop\scriptstyle#2{\mathop{%
\hbox{$\smash{\rightarrow\joinrel\relbar\joinrel\rightarrow\joinrel\relbar}$}%
}\limits^{\raise5pt\hbox{\kern.3pt${\downarrow}$}\atop\raise4pt\hbox{$\smash{|}$}}}#3%
}
\newcommand{\IntType}[3]{\!%
\setbox0=\hbox{$\scriptstyle #3\atop\scriptstyle #1\;#2$}%
{\scriptstyle#3\atop\scriptstyle%
\hbox to\wd0{\hfil$\scriptstyle #1$\hfil\hfil$\scriptstyle#2$\hfil}}%
\!}
\newcommand{\qe}{\qed\vskip2ex}
\makeatletter
\@addtoreset{equation}{section}

\makeatother

\def\bl{\begin{lem}\sl}
\def\el{\end{lem}}
\def\bt{\begin{thm}\sl}
\def\et{\end{thm}}
\def\bp{\begin{pro}\sl}
\def\ep{\end{pro}}
\def\br{\begin{rem}\sl}
\def\er{\end{rem}}
\def\bc{\begin{cor}\sl}
\def\ec{\end{cor}}
\def\bd{\begin{dfn}\rm}
\def\ed{\end{dfn}}
\def\bn{\begin{note}\rm}
\def\en{\end{note}}
\def\bex{\begin{example}\rm}
\def\eex{\end{example}}
\def\proof{{\it Proof.}}

\begin{document}
\title{\Large The automorphism group of the $\Z_2$-orbifold of\\ the Barnes-Wall lattice vertex operator algebra of central charge $32$}

\author{Hiroki Shimakura\footnote{The author was partially supported by JSPS KAKENHI Grant Numbers 18740001, 20549004, 23540013.} \\
}
\date{\small{\it Research Center for Pure and Applied Mathematics,\\ Graduate School of Information Sciences,
Tohoku University,\\
Aramaki aza Aoba 6-3-09, Aoba-ku Sendai-city, 980-8579, Japan}\\
{\rm e-mail: shimakura@m.tohoku.ac.jp}\\
\vspace{0.5cm}
{\rm 2010} {\it Mathematics Subject Classification}. Primary 17B69; Secondary 20B25.\\
}

\maketitle

\begin{abstract} 
In this article, we prove that the full automorphism group of the $\Z_2$-orbifold of the Barnes-Wall lattice vertex operator algebra of central charge $32$ has the shape $2^{27}.E_6(2)$.
In order to identify the group structure, we introduce a graph structure on the Griess algebra and show that it is a rank $3$ graph associated to $E_6(2)$. 
\end{abstract}

\section*{Introduction}
The automorphism group of a vertex operator algebra (VOA) relates to many areas in mathematics.
For example, the automorphism group of the moonshine VOA is the Monster simple group, and this realization plays important roles in the moonshine phenomena, an amazing relation between group theory, VOA theory and number theory (\cite{Bo2}).
Hence it has been an important problem to determine the automorphism group of a VOA.

One of the standard methods of determining the automorphism group $\Aut V$ of a VOA $V$ is to restrict its action to the non-trivial lowest weight subspace.
If the weight $1$ subspace $V_1$ is non-zero, then $\Aut V$ can be studied by using the Lie algebra structure on $V_1$.
If the weight $1$ subspace is zero, then the weight $2$ subspace $V_2$ forms a commutative (non-associative) algebra, called the Griess algebra.
Indeed, the automorphism groups of some VOAs were calculated by using the Griess algebras (\cite{DG0,DG,FLM,Gr1,MM, LSY}).
However, in general, it is difficult to determine the automorphism group of the Griess algebra.

\medskip

The Barnes-Wall lattices ${\rm BW}_{2^n}$ of rank $2^n$ are even lattices in $\R^{2^n}$ (\cite{BW}). They have remarkable automorphism groups of shape $2^{1+2n}_+.\Omega^+(2n,2)$ if $n\ge4$, and play important roles in lattice theory, such as sphere packings (cf.\ \cite{CS}).
Since there are many analogues between lattice theory and VOA theory, it is believed that VOAs with large symmetry can be obtained from Barnes-Wall lattices; 
indeed, the VOAs $V_{\sqrt2 {\rm BW}_8}^+$ and $V_{\sBW}^+$ are such examples (\cite{Sh2,Sh3,Ho}).
The next example would be the holomorphic VOA $\tilde{V}_{\BW}$ of central charge $32$ obtained as the $\Z_2$-orbifold of the lattice VOA $V_{\BW}$;
in \cite{Mi3,Mi}, an infinite series of holomorphic framed VOAs with finite automorphism groups are constructed, and $\tilde{V}_{\BW}$ turns out to be the VOA next to the moonshine VOA.

The automorphism group of the lattice VOA $V_L$ associated to even lattice $L$ was determined in \cite{DN} as $N(V_L)\cdot O(\hat{L})$, where $N(V_L)=\langle \exp(a_0)\mid a\in (V_L)_1\rangle $ and $O(\hat{L})$ is the automorphism group of the central extension $\hat{L}$ of $L$ by $\Z/2\Z$.
On the other hand, the automorphism group of the $\Z_2$-orbifold $\tilde{V}_L$ has not been well studied except for the moonshine VOA (cf.\ \cite{FLM}).

\medskip

The main theorem of this article is the following:
\setcounter{section}{2}
\setcounter{dfn}{29}
\bt Let $\tilde{V}_{\BW}$ be the $\Z_2$-orbifold of the lattice VOA associated to the Barnes-Wall lattice $\BW$ of rank $32$.
Then its automorphism group has the shape $2^{27}.E_6(2)$.
\et
\setcounter{section}{0}

Note that this theorem was conjectured by Griess in a workshop held on July 2004 at Edinburgh.

\medskip

Let us provide two related topics, VOAs of class $S^n$ and trace functions.

A VOA is said to be of class $S^n$ if its trivial component with respect to the full automorphism group coincides with the subVOA generated by the conformal element up to weight $n$ (\cite{Ma}).
This is a sufficient condition that all homogeneous subspaces are conformal $n$-designs, introduced in \cite{Ho}; the notion of conformal designs based on VOAs is a natural analogue of the notion of block and spherical designs based on binary codes and lattices, respectively.
Although the VOAs of class $S^6$ with $V_1\neq0$ are classified in \cite{Ho}, few candidates of VOAs of class $S^6$ with $V_1=0$ are known (\cite{Ma,Ho}); $\tilde{V}_{\BW}$ is one of them (\cite[Example 3.2]{Ho}).

For a VOA $V=\oplus_{i=0}^\infty V_i$ and $g\in \Aut V$, a trace function is given by $T_g(q)=\sum_{i=0}^\infty {\rm Tr}(g_{|V_i})q^{i-c/24}$.
For example, trace functions for the moonshine VOA provide a correspondence between the conjugacy classes of the Monster simple group and certain modular functions (\cite{Bo2}).
Hence, it would be an interesting problem to study similar phenomena on trace functions for $\tilde{V}_{\BW}$ and its automorphism group of shape $2^{27}.E_6(2)$; notice that $T_1(q)$ is just the graded dimension and $T_1(q)=\sqrt[3]{j}(j-992)$, where $\sqrt[3]{j}=q^{-1/3}(1+248q+4124q^2+\dots)$ (cf.\ \cite{HoPhd}).

\medskip

Let $G$ be the automorphism group of $\tilde{V}_{\BW}$.
The proof of the main theorem is divided into three steps:

\begin{enumerate}[(i)]
\item Construct an elementary abelian $2$-subgroup $A$ of $G$ of order $2^{27}$;
\item Prove that $N_G(A)/A$ is isomorphic to $E_6(2)$, where $N_G(A)$ is the normalizer of $A$ in $G$;
\item Prove that $A$ is normal in $G$.
\end{enumerate}
Let us explain these steps in more detail.

(i) Set $V=\tilde{V}_{\BW}$ and $U=V_{\sBW}^+$.
By the description of $\BW$ as an overlattice of $\sBW\oplus\sBW$ in \cite{Gr3}, $V$ contains $U^{\otimes2}$ as a subVOA such that 
as $U^{\otimes 2}$-modules, $$V\cong\bigoplus_{[M]\in S}M^{\otimes 2},$$ where $S$ is the set of all isomorphism classes of irreducible $U$-modules.
Since $S$ forms an elementary abelian $2$-group of order $2^{10}$ under the fusion rules (\cite{AD,ADL,Sh2}), $V$ is a simple current extension of $U^{\otimes2}$ graded by $S$.
Hence the dual $S^*$ of $S$ acts on $V$ as an automorphism group.
By \cite[Theorem 3.3]{Sh2} the restriction group homomorphism 
$$ \rho:N_G(S^*)\to\{g\in\Aut U^{\otimes 2}\mid \{g\circ (M^{\otimes 2})\mid M\in S\}=\{M^{\otimes2}\mid M\in S\}\}$$
is surjective and its kernel is identical to $S^*$, where  $N_G(S^*)$ is the normalizer of $S^*$ in $G$;
the image of $\rho$ is described by $\Aut U^{\otimes 2}\cong\Aut U\wr\Z_2$ (\cite{Sh6}) and the description of the action of $\Aut U$ on $S$ in \cite{Sh2}.
Let us consider the subgroup $A$ of $G$ generated by $S^*$ and $\rho^{-1}(\langle \Delta (O_2(\Aut U)),\tau_U\rangle)$, where $\Delta(O_2(\Aut U))$ is the subgroup $\{(h,h)\mid h\in O_2(\Aut U)\}$ of $\Aut U\wr\Z_2$ and $\tau_U\in\Aut U^{\otimes2}$ sends $x\otimes y$ to $y\otimes x$; one can see that $A$ is abelian, of order $2^{27}$ and normal in $N_G(S^*)$.
We view $A$ as a subgroup of $C_G(z)$, 
where $z$ is the involution in $G$ associated to the $\Z_2$-grading $V=V_{\BW}^+\oplus V_{\BW}^{T,+}$.
By using the explicit description of $C_G(z)$ in terms of $\Aut V_{\BW}^+$ (\cite{FLM}), we prove that $A$ is normal in $C_G(z)$ and that $A$ is elementary abelian.

(ii) By the arguments in (i), $N_G(A)$ contains both $C_G(z)$ and $N_G(S^*)$.
In order to identify the structure of $N_G(A)$, we consider its action on the graph $\Gamma$ defined as follows: the vertex set $V\Gamma$ is the set of eigenspaces of $A$ in the orthogonal complement $B$ of the conformal element in $V_2$, and two distinct vertices $B_1$ and $B_2$ are adjacent if there exist $x_i\in B_i$ such that $x_1\cdot x_2\neq0$, where $\cdot$ is the product in $V_2$;
$N_G(A)/A$ acts on this graph as an automorphism group.
By the action of $N_G(S^*)$ and $C_G(z)$ on $B$, $B$ is an irreducible module for $N_G(A)$.
Hence $N_G(A)/A$ acts transitively on $V\Gamma$.
For the irreducible character $\mu\in A^*$ specified in (\ref{Def:chi}), we show that the dimension of the eigenspace $B_\mu$ with character $\mu$ is $1$ and its stabilizer in $N_G(A)/A$ is identical to $N_G(S^*)/A$.
Moreover, we determine the vertices adjacent to $B_\mu$.
Hence, the cardinality of $V\Gamma$ is equal to $\dim B(=139503)$, the valency is $4590$, and the stabilizer of a vertex in $N_G(A)/A$ is isomorphic to $N_G(S^*)/A$.
In addition, we also prove in Proposition \ref{LG3} that the action of $N_G(A)/A$ on $V\Gamma$ is rank $3$.
According to the classification of primitive permutation groups of rank $3$ (\cite{Li} and references therein), $N_G(A)/A$ must be $E_6(2)$.

(iii) Since $V$ is finitely generated, $G$ is algebraic (\cite{DG2}).
By an argument similar to \cite{Ti}, $G$ is discrete, and hence $G$ is finite.
Using finite group theory, we prove that $A$ is normal in $G$ (Proposition \ref{PNormal});
remark that the classification of finite simple groups is used to prove the following: there are no finite simple groups $K$ with an involution $t$ such that $C_K(t)\cong 2^{1+32}.2^{10}.\Omega^+(10,2)$.

\medskip

The organization of this article is as follows:
In \S 1, we recall and give some definitions and facts necessary in this article, including a review of $V_{\sBW}^+$ in \S 1.5.
In \S 2.1 and \S 2.2, we discuss the structure of two subgroups $C_G(z)$ and $N_G(S^*)$ of $G$ and their actions on $B$, respectively.
In \S 2.3, we construct an elementary abelian $2$-subgroup $A$ of order $2^{27}$, which is a normal subgroup of both $C_G(z)$ and $N_G(S^*)$.
In \S 2.4, we introduce a graph $\Gamma$, and in \S 2.5, we study its structure and the action of $N_G(A)/A$ on it.
In \S 2.6, we prove that the action of $N_G(A)/A$ on $\Gamma$ is rank $3$, and in \S 2.7 we review the classification of permutation groups of rank $3$ and identify $N_G(A)/A$.
In \S 2.8, we prove that $A$ is normal in $G$.

\paragraph{Notations}
\begin{small}
\begin{center}
\begin{longtable}{ll}
$\1_V$& the vacuum vector of a VOA $V$.\\
$2^{n}$& an elementary abelian $2$-group of order $2^n$.\\
$2^{1+2m}_+$& an extraspecial $2$-group of order $2^{1+2m}$ of plus type.\\
$A$& the elementary abelian $2$-subgroup of $\Aut V$ defined in (\ref{ac}).\\
$A^*$& the set of all irreducible characters of $A$.\\
$\Aut X$& the automorphism group of $X$, where $X$ is a lattice, a graph, or a VOA.\\
${\rm BW}_{2^n}$& the Barnes-Wall lattice of rank $2^n$.\\
${{B}}$& the orthogonal complement of the conformal element $\w_V$ in the Griess algebra $V_2$.\\
${{B}}_\chi$& the eigenspace of $A$ in ${{B}}$ with character $\chi$.\\
$C$& the subgroup of $\Aut U$ which acts trivially on $S$, and in fact, $C=O_2(\Aut U)$.\\
$C_H(K)$& the centralizer of $K$ in a group $H$.\\
$\Delta(H)$& the subgroup $\{(h,h)\in H\times H\mid h\in H\}$ of $H\wr\Z_2$.\\
$G$& the automorphism group of $V=\tilde{V}_{\BW}$ in Section 2.\\
$\Gamma$& the graph associated to ${B}$ and ${A}$ (see Definition \ref{Def1}).\\
$\eta$& the canonical group homomorphism from $O(\hat{L})$ to $\Aut L/\langle-1\rangle$ defined in (\ref{Seq2}).\\
$L^*$& the dual lattice of $L$.\\
$L(n)$& the subset $\{a\in L\mid \langle a,a\rangle=n\}$ of a lattice $L$.\\
$[M]$& the isomorphism class of a module $M$ for a VOA.\\
$N_H(K)$& the normalizer of a subgroup $K$ in a group $H$.\\
$O_2(H)$& the maximal normal $2$-subgroup of a group $H$.\\
$O(\hat{L})$& the automorphism group of the central extension $\hat{L}$ of $L$ by $\Z/2\Z$\\
$\rho$& the restriction group homomorphism from $N_G(S^*)$ to $\Aut U^{\otimes2}$ defined in (\ref{ng}).\\
$S=R(U)$& the set of all isomorphism classes of irreducible $U$-modules for the VOA $U$.\\
$S^*$& the subgroup of $\Aut V$ induced from the dual group of $S$.\\
$\theta$& an involution in $O(\hat{L})$ whose image under the map $O(\hat{L})\to \Aut L$ is $-1$\\ 
$\omega_V$& the conformal element of a VOA $V$.\\
$U$, $V$& $U=V_{\sBW}^+$ and  $V=\tilde{V}_{\BW}$ in Section 2.\\
$V\Gamma$& the set of all vertices in a graph $\Gamma$.\\
$V_L^+$& the subVOA of the lattice VOA $V_L$ fixed by $\theta$.\\
$\tilde{V}_L$& the $\Z_2$-orbifold of the holomorphic lattice VOA $V_L$ associated to $\theta$.\\
$\nu$& the canonical group homomorphism from $C_{\Aut \tilde{V}_L}(z)$ to $\Aut V_L^+$ defined in (\ref{Seq1}).\\
$\mu$& the irreducible character of $A$ defined in (\ref{Def:chi}).\\
$X.Y$& a group extension with normal subgroup $X$ and quotient $Y$.\\
$z$& the automorphism of $\tilde{V}_L$ which acts by $1$ and $-1$ on $V_L^+$ and $V_L^{T,+}$, respectively.\\
$\Z_2$& a cyclic group of order $2$.\\
$Z(H)$& the center of a group $H$.\\
$\Omega^+(10,2)$& the simple normal subgroup of the orthogonal group $O^+(10,2)$ of degree $10$ over $\F_2$.
\end{longtable}
\end{center}
\end{small}

\section{Preliminaries}
In this section, we recall or give some definitions and facts necessary in this article.

\subsection{Barnes-Wall lattices}
In this subsection, we review the Barnes-Wall lattices of rank $16$ and $32$ (cf.\ \cite{BW,CS,Gr3}).

Let us start by recalling an explicit construction of the Barnes-Wall lattice $\sBW$ of rank $16$.
Let $\alpha_1,\alpha_2,\dots,\alpha_{16}$ be an orthogonal basis of $\R^{16}$ of (squared) norm $2$.
Let $\RM(1,4)$ denote the first order Reed-Muller code of length $16$, that is, $\RM(1,4)$ is a $5$-dimensional subspace of $\F_2^{16}$ spanned by $(1^{16})$, $((10)^8)$, $((1100)^4)$, $(1^40^41^40^4)$, and $(1^80^8)$.
For an element $c=(c_i)\in\F_2^{16}$, set $\alpha_c=\sum_{i=1}^{16} c_i\alpha_i$, where we view $c_i$ as an element in $\{0,1\}$.
Then 
\begin{equation}
\sBW\cong\sum_{1\le i,j\le 16}\Z(\alpha_i+\alpha_j)+\sum_{c\in \RM(1,4)}\Z\frac{\alpha_c}{2}.\label{BW16}
\end{equation}
Notice that $\sBW^*/\sBW\cong 2^8$ and $\Aut \sBW\cong 2^{1+8}_+.\Omega^+(8,2)$ (cf.\ \cite[\S 10 in Chapter 4]{CS}).

Let $\BW$ denote the Barnes-Wall lattice of rank $32$.
It is well known that $\BW$ is even and unimodular, it is generated by vectors of norm $4$ and $\Aut \BW$ is of shape $2^{1+10}_+.\Omega^+(10,2)$ (cf.\ \cite{Gr3} and references therein).
In \cite[Definition 6.2]{Gr3}, $\BW$ was constructed from $\sBW$  as follows:
\begin{equation}
\BW\cong\sBW\oplus\sBW+\{(v,v)\in\R^{32}\mid v\in \sBW^*\},\label{BW32}
\end{equation}
where $\sBW^*$ is the dual lattice of $\sBW$.
Under the isomorphism (\ref{BW32}), the sublattice $\BW[1]$ of $\BW$ was defined in \cite[Definition 6.2 and Lemma 7.6]{Gr3} by 
\begin{equation}
\BW[1]\cong 2\sBW^*\oplus2\sBW^*+\{(v,v)\in\R^{32}\mid v\in\sBW\}.\label{Def:K}
\end{equation}
It follows from $\sqrt2\sBW^*\cong\sBW$ that $\BW[1]\cong\sqrt2\BW$.
Hence $\BW/\BW[1]\cong 2^{16}$.
By \cite[Corollary 6.5]{Gr3}, $\Aut\BW$ acts on $\BW[1]$ as $\Aut \BW[1]$.
By the same procedure,  we obtain the sublattice $(\BW[1])[1]$ of $\BW[1]$; clearly $(\BW[1])[1]=2\BW$.

\bl\label{LBW32} {\rm (cf.\ \cite[Corollary 6.5 and Lemmas 7.6 and 7.13]{Gr3})}
The automorphism group of $\BW$ preserves $\BW[1]$, and the subgroup of $\Aut\BW$ which acts trivially on $\BW[1]/2\BW$ is $O_2(\Aut \BW)\cong 2^{1+10}_+$.
Moreover, $\BW/\BW[1]$ and $\BW[1]/2\BW$ are irreducible modules over $\F_2$ for $\Aut\BW$.
\el

The following lemma can be proved by the explicit construction of irreducible modules for extraspecial $2$-groups of plus type (cf.\ \cite[\S 5.5]{FLM}).

\bl\label{L12} Let $H\cong 2^{1+2m}_+$ and let $W$ be a $2^m$-dimensional irreducible $H$-module over $\C$ on which the central element of $H$ acts by $-1$.
Then $W\otimes_\C W$ is a $2^{2m}$-dimensional irreducible $H*H$-module by the natural action $(g,h):v\otimes w\mapsto g(v)\otimes h(w)$, where $*$ is the central product and $H*H\cong 2^{1+4m}_+$.
\el

Later, we will apply this lemma to the $2^4$-dimensional irreducible module $\C\otimes_\Z\sBW$ for $O_2(\Aut\sBW)\cong 2^{1+8}_+$.

\subsection{Vertex operator algebras and modules}
In this subsection, we recall the notion of vertex operator algebras (VOAs) and modules from \cite{Bo,FLM,FHL}.
Throughout this article, the ground field is $\C$ if not otherwise stated.

A vertex operator algebra (VOA) $(V,Y,\1_V,\omega_V)$ is a $\Z_{\ge0}$-graded
 vector space $V=\bigoplus_{m\in\Z_{\ge0}}V_m$ equipped with a linear map

$$Y(a,z)=\sum_{i\in\Z}a_{(i)}z^{-i-1}\in ({\rm End}\, V)[[z,z^{-1}]],\quad a\in V,$$
the {\it vacuum vector} $\1_V$ and the {\it conformal element} $\omega_V$ satisfying a number of conditions (\cite{Bo,FLM}).
We often denote it by $V$ or $(V,Y)$ simply.

Two VOAs $(V,Y,\1_V,\omega_V)$ and $(V^\prime,Y',\1_{V'},\omega_{V'})$ are said to be {\it isomorphic} if there exists a linear isomorphism $g$ from $V$ to $V^\prime$ satisfying $$ g\omega_V=\omega_{V'}\quad {\rm and}\quad gY(v,z)=Y'(gv,z)g\quad {\rm for}\ \forall v\in V.$$ 
When $V=V'$, such an isomorphism is called an {\it automorphism}.
The group of all automorphisms of $V$ is called the {\it automorphism group} of $V$ and is denoted by $\Aut V$.

A {\it vertex operator subalgebra} (or a {\it subVOA}) is a graded subspace of
$V$ which has a structure of a VOA such that the operations and its grading
agree with the restriction of those of $V$ and that they share the vacuum vector.
When they also share the conformal element $\omega_V$, we will call it a {\it full subVOA}.

An (ordinary) module $(M,Y_M)$ for a VOA $V$ is a $\C$-graded vector space $M=\bigoplus_{m\in\C} M_{m}$ equipped with a linear map
$$Y_M(a,z)=\sum_{i\in\Z}a_{(i)}z^{-i-1}\in ({\rm End\ }M)[[z,z^{-1}]],\quad a\in V$$
satisfying a number of conditions (\cite{FHL}).
We often denote it by $M$ and its isomorphism class by $[M]$.
The {\it weight} of a homogeneous vector $v\in M_k$ is $k$.
If $M$ is irreducible then there exists the {\it lowest weight} $h\in\C$ of $M$ such that $M=\bigoplus_{m\in\Z_{\ge0}}M_{h+m}$ and $M_h\neq0$.

A VOA is said to be {\it simple} if $V$ is irreducible as a $V$-module.
A VOA is said to be  {\it rational} if any (admissible) module is completely reducible.
A rational VOA is said to be {\it holomorphic} if itself is the only irreducible module up
to isomorphism.
A VOA is said to be {\it of CFT type} if $V_0=\C\1_V$, and is said to be {\it $C_2$-cofinite} if $\dim V/{\rm Span}_\C\{u_{(-2)}v\mid u,v\in V\}<\infty$.
Let $R(V)$ denote the set of all isomorphism classes of irreducible $V$-modules.

Let $M_a$ and $M_b$ be modules for a simple, rational and $C_2$-cofinite VOA $V$ of CFT type.
Then, up to isomorphism, there exists a module $M_a\boxtimes_V M_b$ for $V$ satisfying some conditions, called the {\it fusion product} (cf.\ \cite{HL}).
Note that $M_a\boxtimes_V M_b$ is uniquely determined by the isomorphism classes of $M_a$ and $M_b$, up to isomorphism, and that its $V$-module structure is described by the fusion rules (cf.\ \cite{FHL}).

Let $M=(M,Y_M)$ be a module for a VOA $V$.
For $g\in\Aut V$, let $g\circ M=(M,Y_{g\circ M})$ denote the $V$-module defined by $Y_{g\circ M}(v,z)=Y_M(g^{-1}v,z)$.
Note that if $M$ is irreducible then so is $g\circ M$ and that if $M_a\cong M_b$ then $g\circ M_a\cong g\circ M_b$.
Hence $\Aut V$ acts on $R(V)$.
Since this action preserves the fusion rules (cf.\ \cite[Lemma 1.3]{Sh6}), it also does the fusion product.

Let $V=\bigoplus_{i=0}^\infty V_i$ be a VOA such that $V_0=\C\1_V$ and $V_1=0$.
Consider the first non-trivial subspace $V_2$ and set
$$ a\cdot b=a_{(1)}b,\quad (a,b)\1_V=a_{(3)}b,\quad {\rm for}\ a,b\in V_2.$$
Then $\cdot$ is a commutative (non-associative) product and $(,)$ is an invariant bilinear form.
The space $V_2$ equipped with them is called the {\it Griess algebra} of $V$.
Since $\Aut V$ fixes $\omega_V$ and preserves $(,)$, it also acts on $\omega_V^\perp=\{v\in V_2\mid (v,\omega_V)=0\}$.

\subsection{Simple currents and simple current extensions}
In this subsection, we recall the notion of simple currents and simple current extensions.

Let $V(0)$ be a simple, rational and $C_2$-cofinite VOA of CFT type.
An irreducible $V(0)$-module $M_a$ is called a {\it simple current} if for any irreducible $V(0)$-module $M_b$, there exists a unique irreducible $V(0)$-module $M_c$ satisfying $M_a\boxtimes_{V(0)} M_b\cong M_c$, up to isomorphism.
Let $\{V(\alpha)\mid \alpha\in D\}$ be a set of inequivalent irreducible $V(0)$-modules indexed by an abelian group $D$.
A simple VOA $V=\bigoplus_{\alpha\in D}V(\alpha)$ is called a {\it simple current extension} of $V(0)$ if it carries a $D$-grading and every $V(\alpha)$ is a simple current.
Note that $V(\alpha)\boxtimes_{V(0)}V(\beta)\cong V({\alpha+\beta})$.

\begin{pro}{\rm (\cite[Proposition 5.3]{DM})}\label{PSCE} Let $V(0)$ be a simple, rational and $C_2$-cofinite VOA of CFT type and let $V=\bigoplus_{\alpha\in D}V(\alpha)$ and $\tilde{V}=\bigoplus_{\alpha\in D}\tilde{V}(\alpha)$ be simple current extensions of $V(0)=\tilde{V}(0)$.
If $V(\alpha)\cong \tilde{V}(\alpha)$ as $V(0)$-modules for all $\alpha\in D$, then $V$ and $\tilde{V}$ are isomorphic.
\end{pro}

By the $D$-grading of a simple current extension $V=\bigoplus_{\alpha\in D}V(\alpha)$ of $V(0)$, the group $D^*$ consisting of irreducible characters of $D$ acts on $V$ as an automorphism group; $\chi\in D^*$ acts on $V(\alpha)$ by the scalar multiplication $\chi(\alpha)$ for each $\alpha\in D$.
The normalizer of $D^*$ in $\Aut V$ is describe as follows:

\bp\label{PSh2} {\rm (\cite[Theorem 3.3]{Sh2})} The restriction homomorphism 
\begin{equation}
N_{\Aut V}(D^*)\to\biggl\{ g\in\Aut V(0)\ \biggl|\  \{g\circ[V(\alpha)]\mid \alpha\in D\}=\{[V(\alpha)]\mid \alpha\in D\}\biggr\}\label{eta}
\end{equation}
is surjective, and its kernel is $D^*$.
\ep

\bl\label{LD*} Let $V$ be a simple current extension of $V(0)$.
Then $N_{\Aut V}(D^*)=\{g\in\Aut V\mid g(V(0))=V(0)\}$.
\el
\proof\ Since $V(0)$ is the subVOA of $V$ fixed by $D^*$, $N_{\Aut V}(D^*)$ preserves $V(0)$.
Let $g\in \Aut V$ with $g(V(0))=V(0)$.
Since $V$ is a simple current extension of $V(0)$, the multiplicity of any irreducible $V(0)$-submodule of $V$ is $1$.
Hence $g$ permutes irreducible $V(0)$-submodules of $V$, and thus $g\in N_{\Aut V}(D^*)$.
\qe

\subsection{Lattice VOA and the $\Z_2$-orbifold}
Let $L$ be an even lattice of rank $n$ and let $V_L$ denote the lattice VOA associated with $L$ (\cite{Bo,FLM}).
Let $\hat{L}$ be the central extension of $L$ by $\{\pm1\}\cong\Z/2\Z$ with commutator relation $[e^\alpha,e^\beta]=(-1)^{\langle \alpha,\beta\rangle}$.
The automorphism group $O(\hat{L})$ of $\hat{L}$ is described as follows (\cite[Proposition 5.4.1]{FLM})
\begin{equation}
1\to \Hom(L,\Z_2)\to O(\hat{L}) \overset{\gamma}{\to} \Aut L\to 1;
\end{equation}
notice that $O(\hat{L})$ acts on $V_L$ as an automorphism group (\cite{FLM,DN}).
We often call $g\in O(\hat{L})$ a {\it lift} of $h\in\Aut L$ if $\gamma(g)=h$.
Let $\theta\in \gamma^{-1}(-1)$, a lift of $-1\in\Aut L$.
Let $V_L^+$ denote the subspace of $V_L$ fixed by $\theta$; it is a simple VOA of CFT type of central charge $n$; notice that the VOA structure of $V_L^+$ is independent of the choice of $\theta$ up to isomorphism (\cite[Corollary D.7]{DGH}).
Clearly, $O(\hat{L})/\langle\theta\rangle\subset\Aut V_L^+$ and 
\begin{equation}
1\to \Hom(L,\Z_2)\to O(\hat{L})/\langle\theta\rangle \overset{\eta}{\to} \Aut L/\langle-1\rangle\to 1\label{Seq2}
\end{equation}
is exact.
Indeed $\Hom(L,\Z_2)$ is a subgroup of $\Aut V_L$ which commute with $\theta$.
The explicit action of $\Hom(L,\Z_2)$ on $V_L$ shows the following lemma (cf.\ \cite[\S 10.4]{FLM}):

\bl\label{Lfixpt} Let $f\in\Hom(L,\Z_2)\subset \Aut V_L$.
Then the following hold:
\begin{enumerate}[\rm (1)]
\item $f$ acts trivially on the subspace $\{h(-1)\1_V\mid h\in\C\otimes_\Z L\}(\cong \C\otimes_\Z L)$ of $(V_L)_1$;
\item The subspace of $V_L^+$ fixed by $f$ is $V_N^+$, where $N=\{v\in L\mid f(v)=0\}$.
\end{enumerate}
\el
We often identify $\Hom(L,\Z_2)$ with $L^*/2L^*$ as follows: 
\begin{equation}
L^*/2L^*\to \Hom(L,\Z_2),\quad v+2L^*\mapsto \langle v,\cdot\rangle,\label{Eq:dual}
\end{equation}
where $\langle,\rangle$ is the inner product on $\R\otimes_{\Z} L$.

Assume that $L$ is unimodular, that is, the dual lattice $L^*$ of $L$ is identical to $L$.
Then $V_L$ is holomorphic (\cite{Do}).
By \cite[Theorem 1.2]{DLM}, $V_L$ possesses a unique irreducible $\theta$-twisted module $V_L^T$, up to isomorphism.
It follows from the explicit construction of an irreducible $\theta$-twisted module in \cite[Chapter 9]{FLM} that the weights of $V_L^T$ are half-integers.
Let $V_L^{T,+}$ be the subspace of $V_L^T$ with integral weights; it is an irreducible module for $V_L^+$.
Set $$ \tilde{V}_L = V_L^+ \oplus V_L^{T,+}.$$ It is well known that $ \tilde{V}_L$
has a unique holomorphic VOA structure by extending its
$V_L^+$-module structure (cf.\ \cite{FLM,DGM}).
Note that the central charge of $\tilde{V}_L$ is $n$.
The VOA  $\tilde{V}_L$ is often called the {\it $\Z_2$-orbifold}\footnote{In some literature, $V_L^+$ is called the $\Z_2$-orbifold of $V_L$ associated to $\theta$.} of $V_L$ associated to $\theta$.
It is known that $\tilde{V}_L$ is simple and of CFT type and that it is a simple current extension of $V_L^+$ graded by $\Z_2$ (cf.\ \cite{ADL}).
Note that if $n\ge 24$ and the minimum norm of $L$ is greater than or equal to $4$, then $(\tilde{V}_L)_1=0$.

Let $z$ denote the automorphism of $\tilde{V}_L$ which acts by $1$ and $-1$ on $V_{L}^+$ and $V_{L}^{T,+}$ respectively.
Since $V$ is a simple current extension of $V_{L}^+$, we have $z\in\Aut \tilde{V}_L$.
By Proposition \ref{PSh2} we obtain the exact sequence 
\begin{eqnarray}
1\to \langle z\rangle\to C_{\Aut\tilde{V}_L}(z)\ \overset{\nu}{\to}\ \Aut V_{L}^+\to 1.\label{Seq1}
\end{eqnarray}
Notice that $\Aut V_L^+=O(\hat{L})/\langle\theta\rangle$ (\cite[Theorem 4.1]{Sh2}) and $\nu^{-1}(\Hom(L,\Z_2))\cong 2^{1+n}_+$ (cf.\ \cite[\S 10.4]{FLM}) for even unimodular lattice $L$ of rank $n$.

Assume that $L(2)=\emptyset$ and $L=\langle L(4)\rangle_{\Z}$; by the same argument as in \cite[Proposition 12.2.5]{FLM}, we can show that $V_{L}^+$ is generated by the weight $2$ subspaces as a VOA.
We also assume that the rank of $L$ is $24$ or $32$; the lowest weight of $V_{L}^{T,+}$ is $2$.
The irreducibility of $V_L^{T,+}$ shows the following:

\begin{lem}\label{Lwt2} Let $L$ be an even unimodular lattice of rank $24$ or $32$.
Assume $L(2)=\emptyset$ and $L=\langle L(4)\rangle_{\Z}$.
Then the VOA $\tilde{V}_L$ is generated by the weight $2$ subspace $(\tilde{V}_L)_2$ as a VOA.
\end{lem}

\subsection{VOA associated to the Barnes-Wall lattice of rank $16$}
In this subsection, we review properties of the VOA $U=V_{\sBW}^+$.

By (\ref{BW16}) $\sBW$ has an orthogonal basis of norm $4$; $U$ is a framed VOA.
By \cite{DGH}, $U$ is rational and $C_2$-cofinite.
Since $\sBW$ is generated by norm $4$ vectors, $U$ is also generated by the weight $2$ subspace of $U$ as a VOA.
It follows from $\sBW(2)=\emptyset$ that $U_1=0$.

Let $S$ denote the set $R(U)$ of all isomorphism classes of irreducible $U$-modules.
By the determination of $R(V_L^+)$ in \cite{AD}, we have $|S|=2^{10}$.
Moreover, by the fusion rules of $V_L^+$ (\cite{ADL}), any irreducible $U$-module is a simple current, and $S$ becomes an elementary abelian $2$-group of order $2^{10}$ under the fusion product.
We view $S$ as a $10$-dimensional vector space over $\F_2$.

Let $q_{U}:S\to \F_2$ be the map defined by setting $q_{U}([M])=1$ if $M$ is $(1/2+\Z)$-graded, and $q_{U}([M])=0$ if $M$ is $\Z$-graded; notice that for $[M]\in S$, the lowest weight of $M$ is $0$, $1$, $3/2$ if $[M]$ is zero, non-zero singular, non-singular, respectively.
It was shown in \cite[Theorem 3.8]{Sh2} that $q_{U}$ is an $\Aut U$-invariant non-singular quadratic form on $S$ of plus type; there exists a natural group homomorphism from $\Aut U$ to the orthogonal group $O(S,q_{U})$; its kernel $C$ is an elementary abelian $2$-group of order $2^{16}$, and its image is the simple subgroup $\Omega^+(10,2)$ of $O(S,q_{U})$ of index $2$ (\cite[Theorem 4.8]{Sh2}); $\Aut U\cong 2^{16}.\Omega^+(10,2)$ and $C=O_2(\Aut U)$.
Since $C$ acts trivially on $S$ and $O(\widehat{\sBW})/\langle\theta\rangle$ is the stabilizer of $[V_L^-]$ (\cite[Propositions 2.8 and 3.10]{Sh2}), $C$ is a subgroup of $O(\widehat{\sBW})/\langle\theta\rangle$.
By \cite[\S 4.3]{Sh2}, the image $\eta(C)$ and kernel of $\eta_{|C}$ in (\ref{Seq2}) are described as follows:
\begin{align}
C\cap \Ker\eta&=\{f\in\Hom(\sBW,\Z_2)\mid f(2\sBW^*)=0\}\cong2^8,\label{Eq:C1}\\
\eta(C)&=O_2(\Aut\sBW)/\langle-1\rangle\cong 2^8.\label{Eq:C2}
\end{align}
By the same argument as in \cite[Proposition 4.16]{Sh6}, we have $\Aut U^{\otimes 2}\cong \Aut U\wr\Z_2$.

Let us consider the action of $\Aut U$ on $\omega_U^\perp$.
The orthogonal complement $\omega_U^\perp$ of the conformal element $\omega_U$ in $U_2$ is a direct sum of $S^2(\C\otimes_\Z \sBW)\cap\omega_U^\perp$ and the subspace $\sum_{\alpha\in \sBW(4)}\C (e^{\alpha}+\theta(e^{\alpha}))$ of the group algebra $\C[\sBW]$ as a vector space, where $S^2(\cdot)$ is the symmetric tensor and $\sBW(4)=\{v\in\sBW\mid \langle v,v\rangle=4\}$.
By $$\dim S^2(\C\otimes_\Z \sBW)\cap\omega_U^\perp=16+\binom{16}{2}-1=135,\ \dim \sum_{\alpha\in \sBW(4)}\C (e^{\alpha}+\theta(e^{\alpha}))=\frac{|\sBW(4)|}{2}=2160,$$
we have $\dim \omega_U^\perp=2295.$
By a similar argument as in \cite[Remark 10.4.13]{FLM}, these two subspaces are irreducible for $O(\widehat{\sBW})/\langle\theta\rangle$.
By \cite[Chapter 11]{FLM} (cf.\ \cite[\S 2.3]{Sh2}) $U$ has an extra automorphism not preserving them.
Hence $\omega_U^\perp$ is irreducible for $\Aut U$.

Let $(\omega_U^\perp)^C=\{v\in \omega_U^\perp \mid g(v)=v\ \text{for\ all}\ g\in C\}$. 
Since $C$ is normal in $\Aut U$, $(\omega_U^\perp)^C$ is a submodule of $\omega_U^\perp$ for $\Aut U$.
Since the action of $C$ on $\omega_U^\perp$ is not trivial on $\omega_U^\perp$ and $\omega_U^\perp$ is irreducible for $\Aut U$, we have $(\omega_U^\perp)^C=0$, that is, the subspace of $U_2$ fixed by $C$ is $\C\omega_U$.

{\it Ising vectors}, vectors in the weight $2$ subspace of a VOA generating a simple Virasoro VOA $L(1/2,0)$ of central charge $1/2$ as a subVOA, are classified for $V_L^+$ with $L(2)=\emptyset$ in \cite{Sh10} as follows: any Ising vector has a form $\omega^\pm(\alpha)$ or $\omega(E,\varphi)$, where $\alpha\in L(4)$, $E$ is a sublattice of $L$ isomorphic to $\sqrt2E_8$ and $\varphi\in\Hom(E,\Z/2\Z)$.
For an Ising vector $e$, the automorphism $\tau_e$ of a VOA , called the $\tau$-involution, is introduced in \cite{Mi0}.

Let $I$ be the set of all Ising vectors in $U$.
Let us show that the normal subgroup $N=\langle\tau_e\mid e\in I\rangle$ is identical to $O_2(\Aut U)$.
Let $I_1$ (resp. $I_2$) be the set of all Ising vectors in $U$ of form $\omega^\pm(\alpha)$ (resp. $\omega(E,\varphi)$); by the classification in \cite{Sh10} $I=I_1\cup I_2$.
By $\Aut U\cong2^{16}.\Omega^+(10,2)$ and $N\neq\{1\}$, $N$ must be $O_2(\Aut U)$ or $\Aut U$.
One can easily see that if $e=\omega^\pm(\alpha)$, then $\tau_e\in C\cap\Ker\eta\subset O_2(\Aut U)$ (cf. \eqref{Eq:C1}).
Hence it suffices to show that $\Aut U$ is transitive on $I$.

The transitivity of $O(\widehat{\sBW})/\langle\theta\rangle$ on $I_1$ follows from that of $\Aut\sBW$ on the set of norm $4$ vectors in $\sBW$.
Let us consider the action of $O(\widehat{\sBW})/\langle\theta\rangle$ on $I_2$.
Let $E$ be a sublattice of $\sBW$ isomorphic to $\sqrt2E_8$ and $\varphi\in\Hom(E,\Z_2)$.
Since a basis of $E$ is contained in one of $\sBW$, there exists $\tilde{\varphi}\in\Hom(\sBW,\Z_2)$ such that its restriction to $E$ is identical to $\varphi$; clearly $\tilde{\varphi}(\omega(E,0))=\omega(E,\varphi)$.
By \cite[Theorem 12.1]{Gr3}, $\Aut\sBW$ also acts transitively on the set of sublattices of $\sBW$ isomorphic to $\sqrt2E_8$.
Thus $O(\widehat{\sBW})/\langle\theta\rangle$ acts transitively on $I_2$.
In addition, $U$ has an extra automorphism $\sigma$ (\cite[Chapter 11]{FLM}) such that $\sigma(I_1)\cap I_2\neq\emptyset$.
Thus $\Aut U$ acts transitively on $I$.

The following proposition is a summary of this subsection:
\begin{pro}\label{POVE}
\begin{enumerate}[{\rm (1)}]
\item The VOA $U$ is simple, rational, $C_2$-cofinite, and of CFT type;
\item $U_1=0$, and $U$ is generated by $U_2$ as a VOA;
\item $(S,q_U)$ is a non-singular $10$-dimensional quadratic space of plus type over $\F_2$.
Moreover, for $[M]\in S$, the lowest weight of $M$ is $0$, $1$, and $3/2$ if $[M]$ is zero, non-zero singular, and non-singular, respectively;
\item $\Aut U\cong 2^{16}.\Omega^+(10,2)$.
Moreover, $C=O_2(\Aut U)$ satisfies (\ref{Eq:C1}) and (\ref{Eq:C2});
\item $C$ is the subgroup of $\Aut U$ acting on $S$ trivially, and $S$ is irreducible for $\Aut U/C$;
\item $\Aut U^{\otimes2}\cong \Aut U\wr\Z_2$;
\item $\dim\omega_U^\perp=2295$ and $\omega_U^\perp$ is irreducible for $\Aut U$;
\item $U_2^{C}=\{v\in U_2\mid g(v)=v\ {\rm for\ all}\ g\in O_2(\Aut U)\}=\C\omega_U$;
\item $C$ is identical to the subgroup generated by all $\tau$-involutions.
\end{enumerate}
\end{pro}

\section{Automorphism group of the VOA $\tilde{V}_{\BW}$}
Throughout this section, we denote $\tilde{V}_{\BW}$, $V_{\sBW}^+$ and $\Aut \tilde{V}_{\BW}$ by $V$, $U$ and $G$, respectively.
Let $B$ denote the orthogonal complement of $\omega_V$ in $V_2$.

\subsection{Centralizer of the involution associated to the $\Z_2$-grading}\label{SecC}
Let $z$ denote the automorphism of $V$ acting by $1$ and $-1$ on $V_{\BW}^+$ and $V_{\BW}^{T,+}$ respectively.
In this subsection, we study the centralizer of $z$ in $G$, and its action on $B$.

Recall from \S 1.4 that $\Aut V_{\BW}^+=O(\widehat{\BW})/\langle\theta\rangle$ and $\nu^{-1}({\rm Hom}_\Z(\BW,\Z_2))\cong 2^{1+32}_+$.
By $\Aut\BW\cong 2^{1+10}_+.\Omega^+(10,2)$, (\ref{Seq2}) and (\ref{Seq1}), we obtain the following proposition.

\bp\label{PCen} The centralizer $C_{G}(z)$ in $G$ of $z$ has the shape $2^{1+32}_+.(2^{10}.\Omega^+(10,2))$.
\ep

\bn In \cite{Mi}, the shape of $C_G(z)$ was determined by a different method.
\en

By the similar arguments as in \cite[Remark 10.4.13]{FLM}, we obtain the following lemma.
\bl\label{LDec1} Under the action of $C_{G}(z)$, ${{B}}$ breaks into the three irreducible submodules $$ B(1):=S^2(\C\otimes_\Z\BW)\cap B, \qquad B(2):=\sum_{\alpha\in \BW(4)}\C (e^\alpha+\theta(e^\alpha)),\qquad B(3):=T$$ of dimensions 
$527,\ 73440,\ 65536,$ respectively.
\el

\subsection{Normalizer of the elementary abelian $2$-subgroup of order $2^{10}$}
Let $U^{\otimes 2}$ be the full subVOA of $V$ based on the embedding of $\sBW\oplus\sBW$ into $\BW$ in (\ref{BW32}), where $U\cong V_{\sBW}^+$.
By (\ref{BW32}) and \cite[Proposition 4.2]{Sh4}, as $U^{\otimes 2}$-modules,
\begin{equation}
V\cong \bigoplus_{[M]\in S} M^{\otimes 2},\label{DecV}
\end{equation}
where $S=R(U)$ is the set of all isomorphism classes of irreducible $U$-modules.
Since any element in $S$ is a simple current, $V$ is a simple current extension of $U^{\otimes 2}$ graded by $S\cong2^{10}$.
Let $S^*$ denote the dual of the group $S$.
Then $S^*$ acts on $V$ as an automorphism group.
In this subsection, we study the structure of the normalizer of $S^*$ in $G$ and its action on $B$.

Set $C=\{g\in\Aut U\mid g\circ U\cong U\ \forall U\in S\}$; by Proposition \ref{POVE} (4) and (5), $C=O_2(\Aut U)\cong 2^{16}$.
Let $\sigma_U$ denote the automorphism of $U^{\otimes 2}$ which sends $x\otimes y$ to $y\otimes x$.
By Proposition \ref{POVE} (5) and (6), the stabilizer of $\{W\otimes W\mid W\in S\}$ in $\Aut U^{\otimes 2}$ is $(C\wr\langle\sigma_U\rangle)\Delta(\Aut U)$, where $\Delta(\Aut U)=\{(h,h)\mid h\in\Aut U\}$ is the subgroup of $\Aut U\times \Aut U$ commuting with $\sigma_U$.
Applying Proposition \ref{PSh2} to (\ref{DecV}), we obtain the following exact sequences:
\begin{eqnarray}
&1\to S^*\to N_{G}(S^*)\overset{\rho}\to(C\wr\langle\sigma_U\rangle)\Delta(\Aut U)\to 1,\label{ng}\\
&1\to S^*\to C_{G}(S^*)\overset{\rho}\to (C\wr\langle\sigma_U\rangle)\to 1,\label{ES:C}
\end{eqnarray}
where $\rho$ is the restriction homomorphism from $N_{G}(S^*)$ to $\Aut U^{\otimes2}$.
Clearly $C\wr\langle\sigma_U\rangle\cap\Delta(\Aut U)=\Delta(C)$.

\bp\label{PNGS*} The normalizer $N_G(S^*)$ and the centralizer $C_G(S^*)$ of $S^*$ in $G$ have the shapes $2^{10}.(2^{16}\wr\Z_2).\Omega^+(10,2)$ and $2^{10}.(2^{16}\wr\Z_2)$, respectively.
\ep

Let us consider the action of $N_G(S^*)$ on ${{B}}$.
By Lemma \ref{LD*}, $N_G(S^*)$ preserves $(U^{\otimes 2})_2$.

\bl\label{L35} Under the action of $N_G(S^*)$, the subspace $(U^{\otimes 2})_2\cap B$ breaks into the two irreducible submodules:
$$ B(4):=\C(\omega_U\otimes\1_U-\1_U\otimes\omega_U),\qquad B(5):=(\omega_U^\perp\otimes\1_U)\oplus(\1_U\otimes\omega_U^\perp)$$
of dimensions $1$ and $4590$, respectively.
\el
\proof\ Since $N_G(S^*)$ preserves $U^{\otimes 2}$ and $\Aut U^{\otimes 2}$ does $B(4)$ by Proposition \ref{POVE} (6), 
$B(4)$ is an (irreducible) submodule for $N_G(S^*)$.
Since $B(5)$ is the orthogonal complement of $B(4)$ in $(U^{\otimes 2})_2\cap B$, it is also a submodule for $N_G(S^*)$.
Hence it suffices to show that $B(5)$ is irreducible.
By Proposition \ref{POVE} (7), $B(5)$ is the direct sum of two irreducible modules $\{v\otimes\1_U+\1_U\otimes v\mid v\in \omega_U^\perp\}$ and $\{v\otimes\1_U-\1_U\otimes v\mid v\in \omega_U^\perp\}$ for $\rho^{-1}(\Delta(\Aut U))$.
Considering the action of $\rho^{-1}(C\times1)$, we have the irreducibility of $B(5)$ for $N_G(S^*)$ by Proposition \ref{POVE} (8).
\qe

\bp\label{PCN} As a $\langle C_G(z),N_{G}(S^*)\rangle$-module, ${{B}}$ is irreducible.
\ep
\proof\ Let $B(i)$ ($1\le i\le 5$) be the subspaces of $B$ given in Lemmas \ref{LDec1} and \ref{L35}.
It suffices to show that $B(1)$, $B(2)$ and $B(3)$ are not preserved by $N_G(S^*)$.

By $(U^{\otimes 2})_2\cap B\subset (V_{\BW}^+)_2\cap B =B(1)\oplus B(2)$ and the dimension of $B(i)$, 
we have $B(5)\cap B(1)\neq 0$ and $B(5)\cap B(2)\neq 0$.
Hence $B(1)$ and $B(2)$ are not preserved by $N_{G}(S^*)$ since $B(5)$ is irreducible.

Let $M$ be an irreducible $U$-module of twisted type with $q_U(M)=0$; $(M\otimes M)_2\cap B(3)\neq0$.
By the transitivity of $\Aut U/C\cong\Omega^+(10,2)$ on the set of non-zero singular vectors in $S$, there exists $g\in N_G(S^*)$ such that $g(M^{\otimes 2})_2\cap (B(1)\oplus B(2))\neq0$.
This shows that $N_G(S^*)$ does not preserve $B(3)$.
\qe

\subsection{Elementary abelian $2$-subgroup of order $2^{27}$}

Based on the exact sequence (\ref{ng}), we consider the following subgroup of $N_G(S^*)$:
\begin{eqnarray}
A=\rho^{-1}(\langle \Delta(C),\sigma_U\rangle).\label{ac}
\end{eqnarray}
It follows from $S^*=\Ker\rho\cong 2^{10}$, $C\cong 2^{16}$ and $|\sigma_U|=2$ that $|A|=2^{27}$.
In this subsection, we prove that $A$ is elementary abelian and is normal in both $N_G(S^*)$ and $C_G(z)$.

Let $\sigma_V$ be an element in $A$ with $\rho(\sigma_V)=\sigma_U$.
Then $\sigma_V^2\in S^*$; the order of $\sigma_V$ is $2$ or $4$.
Later, we will show in Corollary \ref{CA} that the order of $\sigma_V$ is $2$.

\bl\label{Lac} 
\begin{enumerate}
\item The subgroup $\rho^{-1}(\Delta(C))$ of $A$ is an elementary abelian $2$-group;
\item $A$ is abelian;
\item $A$ is normal in $N_{G}(S^*)$.
\end{enumerate}
\el
\proof\ Since $S^*$ acts on each $U^{\otimes 2}$-submodule of $V$ by a scalar and $A\subset C_G(S^*)$, we have $S^*\subset Z(A)$.

Let us show (1) and (2) by describing an explicit section of $\rho_{|A}$ in terms of $\tau$-involutions.
By Proposition \ref{POVE} (9), $C$ is generated by $\tau$-involutions associate to Ising vectors in $U$.
For $(c,1)\in C\times 1\subset \rho(N_G(S^*))$, there exists an Ising vector $e\in U\otimes\1_U\subset V$ such that ${\tau_e}_{|U^{\otimes2}}=(c,1)$.
Set $f=\sigma_V(e)$; $f\in \1_U\otimes U$ and it is orthogonal to $e$ in $V_2$, 
which shows that $\tau_{f}$ commutes with $\tau_e$; $\tau_e\tau_{f}$ is an involution in $G$ and $\rho(\tau_e\tau_f)=(c,c)$; $\rho^{-1}(c,c)$ is elementary abelian.
For $(c',1)\in C\times1$, we obtain $\tau_{e'}$ and $\tau_{f'}$ by the same procedure.
Since $[\tau_e,\tau_{f'}]=[\tau_{e'},\tau_{f}]=1$ and $[\tau_e,\tau_{e'}]=[\tau_f,\tau_{f'}]\in S^*$, we have $[\tau_e\tau_{f}, \tau_{e'}\tau_{f'}]=1$; $\rho^{-1}(\langle(c,c),(c',c')\rangle)$ is elementary abelian.
Moreover, these involutions generate an elementary abelian $2$-subgroup whose image under $\rho$ is $\Delta(C)$.
Hence (1) holds.
It follows from $\sigma_V(\tau_e\tau_f)\sigma_V^{-1}=\tau_f\tau_e=\tau_e\tau_f$ that $\sigma_V\in Z(A)$; (2) holds.

Let us show (3).
Since $A$ contains $S^*$, it suffices to show that $\rho(A)=\Delta(C)\times\langle\sigma_U\rangle$ is normal in $\rho(N_{G}(S^*))=(C\wr\langle\sigma_U\rangle)\Delta(\Aut U)$ (cf.\ (\ref{ng}) and (\ref{ac})). 
Clearly, $\rho(A)$ is normalized by $\Delta(\Aut U)$, and $\Delta(C)$ commutes with $C\wr\langle\sigma_U\rangle$.
Moreover, for $(c,c')\in C\times C$, $(c,c')\sigma_U(c,c')^{-1}=(cc',cc')\sigma_U\in \rho(A)$.
Hence $\rho(A)$ is normal in $\rho(N_{G}(S^*))$.
\qe

Since $A$ contains $z$, it is a subgroup of $C_G(z)$.
In order to show that $A$ is normal in $C_G(z)$, we consider the structure of $A$ in $C_G(z)$.
Recall that $\nu:C_G(z)\to\Aut V_{\BW}^+$ and $\eta:\Aut V_{\BW}^+\to\Aut \BW/\langle-1\rangle$ are defined in (\ref{Seq1}) and (\ref{Seq2}), respectively, and that $\Ker\nu=\langle z\rangle$, $\Ker\eta={\rm Hom}_{\Z}(\BW,\Z_2)\cong 2^{32}$, and $\Ker \eta\circ \nu\cong 2^{1+32}_+$.

\bl\label{LCent1}\begin{enumerate}
\item $\nu(S^*)\cap\Ker\eta=\{f\in\Hom(\BW,\Z_2)\mid f(\sBW\oplus\sBW)=0\}\cong 2^8$.
\item $\nu({A})\cap\Ker\eta=\{f\in\Hom(\BW,\Z_2)\mid f(\BW[1])=0\}\cong2^{16}$, where $\BW[1]$ is the sublattice of $\BW$ given in (\ref{Def:K}).
\item $\eta\circ\nu({A})=O_2(\Aut \BW/\langle-1\rangle)\cong 2^{10}$.
\item $\eta\circ\nu(S^*)=\langle (-1,1),-1\rangle/\langle-1\rangle$, where $(-1,1)$ is the element in $\Aut\BW$ which sends $(x,y)\mapsto (-x,y)$ in the description (\ref{BW32}) of $\BW$.
\end{enumerate}
\el
\proof\ (1) follows from $V^{S^*}=U^{\otimes 2}$ and Lemma \ref{Lfixpt} (2).

Let $L$ be the sublattice of $\BW$ satisfying $\nu(A)\cap\Ker\eta=\{f\in\Hom(\BW,\Z_2)\mid f(L)=0\}$ and $2\BW\subset L$.
Set $H=\nu\circ \rho^{-1}(\Delta(C\cap \Hom(\sBW,\Z_2)))\cap\Ker\eta$.
It follows from (\ref{ac}) and $\rho^{-1}(\Delta(C))\sigma_V\cap\Ker\eta=\emptyset$ that $\nu({A})\cap\Ker\eta=\langle H, \nu({S^*})\cap\Ker\eta\rangle$.
By the definition of $H$ and (\ref{Eq:C1}), for $f\in\Hom(\BW,\Z_2)$, $f\in H$ if and only if $f(v)=0$ for all $v\in 2\sBW^*\oplus 2\sBW^*+\{(v,v)\mid v\in\sBW^*\}$.
By (1) $$L=(\sBW\oplus\sBW)\cap (2\sBW^*\oplus 2\sBW^*+\{(v,v)\mid v\in\sBW^*\})=\BW[1].$$
Hence by (\ref{Def:K}), we obtain $L=\BW[1]$, which proves (2).

By (2) and $|\nu(A)|=2^{26}$, we have $|\eta\circ\nu({A})|=2^{10}$.
By Lemma \ref{Lac} (2), $\nu({A})$ is abelian, and hence it acts trivially on $\nu({A})\cap\Ker\eta$.
By (\ref{Eq:dual}) and (2), $\eta\circ\nu(A)$ acts trivially on $\BW[1]/2\BW$.
By Lemma \ref{LBW32}, $\eta\circ\nu(A)$ is a subgroup of $O_2(\Aut \BW/\langle-1\rangle)\cong 2^{10}$.
Comparing the orders of the groups, we have (3).

By (1), the order of $\eta\circ\nu(S^*)$ is $2$.
On the other hand, $U^{\otimes 2}$ is the subspace of $V_{\sBW\oplus\sBW}^+$ fixed by some lift of $(-1,1)$.
Hence we have (4).
\qe

\bp\label{PsubG3} The subgroup $A$ is normal in $C_{G}(z)$.
\ep
\proof\ Set $Q=\Ker(\eta\circ\nu)\cap A$.
By Lemmas \ref{LBW32} and \ref{LCent1} (2) and  (\ref{Eq:dual}), $Q$ is normal in $C_G(z)$.
Let us show that $C_G(Q)=A$.
Clearly, $A\subset C_G(Q)$.
Since $Q$ is maximal abelian in $\Ker(\eta\circ\nu)$, we have $\Ker(\eta\circ\nu)\cap C_G(Q)=Q=\Ker(\eta\circ\nu)\cap A$.

Since $C_G(Q)$ acts trivially on $\nu(A)\cap\Ker\eta$, 
so does $\eta\circ\nu (C_G(Q))$ on $\BW[1]/2\BW$ by Lemma \ref{LCent1} (2) and (\ref{Eq:dual}).
By Lemmas \ref{LBW32} and \ref{LCent1} (3), $\eta\circ\nu (C_G(Q))\subset \eta\circ\nu (A)$.
Hence $C_G(Q)\subset A$.
Thus $C_G(Q)=A$.
\qe

\bc\label{CA}
\begin{enumerate}[{\rm (1)}]
\item The normalizer $N_G(A)$ contains both $N_G(S^*)$ and $C_G(z)$.
\item Under the action of $N_G(A)$, ${{B}}$ is irreducible.
\item The group $A$ is an elementary abelian $2$-group of order $2^{27}$. 
\end{enumerate}
\ec
\proof\ (1) follows from Lemma \ref{Lac} (3) and Proposition \ref{PsubG3}, and (2) does from (1) and Proposition \ref{PCN}.
Since $\eta\circ\nu(C_G(z))/O_2(\eta\circ\nu(C_G(z)))\cong\Omega^+(10,2)$ acts irreducibility on $\eta\circ \nu(A)\cong 2^{10}$ and $\eta\circ\nu(\sigma_V)\neq1$, we have $\sigma_V^2=1$.
Hence (3) holds by Lemma \ref{Lac} (1).\qe

\bp\label{PAirr} The subgroup $A$ is an irreducible $N_G(A)/A$-module over $\F_2$.
\ep
\proof\ Let $A_0$ be a non-zero submodule of $A$ for $N_G(A)/A$.
By Corollary \ref{CA} (1), both $C_G(z)$ and $N_G(S^*)$ preserve $A_0$.
 
If $\eta\circ \nu(A_0)\neq 0$, then by the action of $C_G(z)$ and Lemma \ref{LCent1} (3), we have $\eta\circ\nu(A_0)=\eta\circ\nu(A)$.
By the action of $N_G(S^*)$, along with Lemma \ref{LCent1} (1) and (4), $(\Ker\eta)\cap \nu(A_0)\neq0$.
By Lemmas \ref{LBW32} and \ref{LCent1} (2), $\nu(A)=\nu(A_0)$.
By $\Ker\nu=\langle z\rangle\subset S^*$, we have $A_0=A$.

If $\eta\circ\nu(A_0)=0$, then by the same argument above, we have $A_0\cap(\Ker\eta\circ\nu)=A\cap(\Ker\eta\circ\nu)$.
By the action of $N_G(S^*)$, along with Lemma \ref{LCent1} (3) and (4), we have $\eta\circ\nu(A_0)=\eta\circ\nu(A)$, and $A_0=A$.
\qe

\subsection{Graph on eigenspaces in the Griess algebra}
In this subsection, we define a graph on the eigenspaces of $A$ in $B=\omega_V^\perp$, which will be used to determine the group structure of $N_G(A)$.

Let $B=\oplus_{\chi\in A^*} B_\chi$ be the eigenspace decomposition of $A$ in $B$.

\bd\label{Def1} For non-zero distinct eigenspaces $B_1,B_2$, we say that $B_1$ is adjacent to $B_2$ if there exists $x_i\in B_i$ ($i=1,2$) such that $x_1\cdot x_2\neq0$, where $\cdot$ is the product in $B$.
\ed

Let $\Gamma$ denote the graph on the vertex set $V\Gamma=\{B_\chi\mid B_\chi\neq0\}$ defined by Definition \ref{Def1}.
Since the product in the Griess algebra is commutative, $\Gamma$ is undirected.
For $g\in N_G(A)$, $g(B_\chi)={B}_{\chi^g}$, where $\chi^g(a)=\chi(g^{-1}ag)$ for $a\in {A}$; $N_G(A)$ acts on $\Gamma$ as an automorphism group.

\bl\label{LG2} Let $I$ be a subset of $V\Gamma$ and let $H$ be a subgroup of $N_{G}({A})$ preserving ${I}$.
Assume that $\bigoplus_{{B}_\chi\in I}{B}_\chi$ is irreducible for $H$.
Then the following hold:
\begin{enumerate}
\item $H$ is transitive on $I$;
\item For ${B}_1$, ${B}_2\in I$, $\dim {B}_1=\dim {B}_2$.
\end{enumerate}
\el
\proof\ Let $I_0$ be an orbit in $I$ under the action of $H$.
Then $\bigoplus_{{B}_\chi\in I_0}{B}_\chi$ is a submodule of ${B}$ for $H$.
By the irreducibility, we have $I_0=I$, and (1) holds.

For $g\in H$ and ${B}_\chi\in V\Gamma$, $\dim {B}_\chi=\dim {B}_{\chi^g}$.
Hence (2) follows from (1).
\qe

\subsection{Graph structure and the stabilizer of a vertex}
In this subsection, we study the structure of the graph $\Gamma$ defined in the previous subsection and the action of $N_G(A)/A$ on it.

Let $\mu$ denote the irreducible character of $A$ defined by
\begin{eqnarray}
\mu(a)=\left\{\begin{array}{cl}
 \mbox{$1$} & \mbox{${\rm if}\ a\in \rho^{-1}(\Delta(C))$},\\
 \mbox{$-1$} & \mbox{${\rm if}\ a\in \rho^{-1}(\Delta(C))\sigma_V$}.
\end{array}
\right. \label{Def:chi}
\end{eqnarray}
Let us describe the eigenspace ${{B}}_{\mu}$ and its stabilizer in $N_{G}(A)$.

\bl\label{LGr6}
The eigenspace ${{B}}_{\mu}$ is spanned by $\omega_U\otimes\1_U-\1_U\otimes\omega_U$, and $\dim B_{\mu}=1$.
\el
\proof\ Since $\rho^{-1}(\Delta(C))$ acts on $U^{\otimes 2}$ as a subgroup of $\Aut U\times\Aut U$,  it fixes both $\omega_U\otimes\1_U$ and $\1_U\otimes\omega_U$.
Moreover, $\sigma_V(\omega_U\otimes\1_U-\1_U\otimes\omega_U)=-(\omega_U\otimes\1_U-\1_U\otimes\omega_U)$.
Hence $\omega_U\otimes\1_U-\1_U\otimes\omega_U$ is an eigenvector of $A$ with character $\mu$.

By $S^*\subset \rho^{-1}(\Delta(C))$, $B_\mu\subset U^{\otimes 2}=\{u\in V\mid g(u)=u\ \forall g\in S^*\}$.
It follows from $U_0=\C\1_U$ and $U_1=0$ that $(U^{\otimes 2})_2=U_2\otimes\1_U\oplus\1_U \otimes U_2$.
By Proposition \ref{POVE} (8), 
$(U^{\otimes2})^{\Delta(C)}=\C(\omega_U\otimes\1_U)\oplus\C (\1_U\otimes\omega_U)$.
Considering the action of $\sigma_V$, we obtain this lemma.
\qe

\bl \label{LGr1} Let ${{B}}_{\chi}\in V\Gamma$.
Then the following hold:
\begin{enumerate}
\item $B_\chi$ is contained in some $U^{\otimes2}$-submodule of $V$, that is, $B_\chi\subset (M^{\otimes 2})_2$ for some $[M]\in S$;
\item Assume $\chi\neq\mu$.
Then ${{B}}_{\chi}$ is adjacent to ${{B}}_{\mu}$ if and only if ${{B}}_{\chi}\subset (U^{\otimes 2})_2\cap B$.
\end{enumerate}
\el
\proof\ By $S^*\subset A$, ${{B}}_{\chi}$ is an eigenspace of $S^*$; (1) holds.

Let us show (2).
First, we assume that ${{B}}_{\chi}\subset (U^{\otimes 2})_2\cap B$.
Let $x\in {{B}}_{\chi}$ be a non-zero vector.
Then $x=v\otimes \1_U+\1_U\otimes w$ for some $v,w\in U_2$ by $(U^{\otimes 2})_2=U_2\otimes\1_U\oplus\1_U\otimes U_2$.
By $$(\omega_U\otimes\1_U-\1_U\otimes\omega_U)\cdot x=2(v\otimes \1_U-\1_U\otimes w)\neq0,$$
${{B}}_{\mu}$ and ${{B}}_{\chi}$ are adjacent by Lemma \ref{LGr6}.

Next, we assume that ${{B}}_{\chi}\subset (M^{\otimes 2})_2\cap B$ and $M\neq U$.
By Proposition \ref{POVE} (3) and $(M^{\otimes2})_2\neq0$, the lowest weight of $M$ is $1$.
Hence any vector in ${{B}}_{\chi}$ has a form $\sum v_i\otimes w_i$, $v_i,w_i\in M_1$.
By $$(\omega_U\otimes\1_U-\1_U\otimes\omega_U)\cdot \sum v_i\otimes w_i=0,$$
${{B}}_{\mu}$ and ${{B}}_{\chi}$ are not adjacent by Lemma \ref{LGr6}.
\qe

\bp\label{Lgraph2}The stabilizer of ${{B}}_{\mu}$ in $N_{G}(A)$ is equal to ${N_{G}(S^*)}$.
\ep
\proof\ By Lemmas \ref{L35} and \ref{LGr6}, ${N_{G}(S^*)}$ is a subgroup of the stabilizer of ${{B}}_{\mu}$ in $N_{G}(A)$.
Let $g\in N_{G}(A)$ with $g({{B}}_{\mu})={{B}}_{\mu}$.
Let $I$ be a subset of all vertices in $V\Gamma$ adjacent to $B_{\mu}$. 
Then $g$ preserves $\bigoplus_{{{B}}_{\chi}\in I}B_{\chi}$ by the definition of $\Gamma$.
By Lemma \ref{LGr1} (2), ${{B}}_{\mu}\oplus  \bigoplus_{{{B}}_{\chi}\in I}{{B}}_{\chi}=(U^{\otimes 2})_2\cap B$.
Hence $g((U^{\otimes 2})_2)=(U^{\otimes 2})_2$, and $g(U^{\otimes 2})=U^{\otimes 2}$ by Proposition \ref{POVE} (2).
It follows from Lemma \ref{LD*} that $g\in N_{G}(S^*)$.
\qe

\bp\label{LGr7} The subgroup of $N_{G}(A)$ acting trivially on $V\Gamma$ is $C_{G}(A)$.
\ep
\proof\ Clearly, $C_{G}(A)$ acts trivially on $V\Gamma$.
Let $g$ be an element in $N_{G}(A)$ acting on $V\Gamma$ trivially.
Then $\chi^g=\chi$ for all $\chi\in A^*$ with ${{B}}_\chi\neq0$.
It follows from Lemma \ref{Lwt2} that $A$ acts faithfully on ${{B}}$.
Hence, $\langle \chi\in\Irr (A)\mid {{B}}_\chi\neq0\rangle=A^*$, and $g$ commutes with $A$.
\qe

\bp\label{CG2} The subgroup $A$ of $G$ is self-centralizing, that is, $C_{G}(A)=A$.
\ep
\proof\ Since $A$ is abelian, $A\subset C_{G}(A)$.
It follows from $S^*\subset A$ that $C_{G}(A)\subset C_{G}(S^*)$.
It follows from $\sigma_V\in A$ and (\ref{ES:C}) that $\rho(C_{G}(A))\subset \Delta(C)\times\langle\sigma_U\rangle$.
Hence $C_{G}(A)\subset A$ by (\ref{ac}).\qe

Recall that $\dim B=139503$ and $\dim U=2296$ (cf.\ Lemmas \ref{LDec1} and \ref{L35}).
Combining Lemmas \ref{LG2}, \ref{LGr6}, and \ref{LGr1} (2), Propositions \ref{Lgraph2}, \ref{LGr7} and  \ref{CG2}, we obtain the following:

\bt\label{CM} The graph $\Gamma$ has $139503$ vertices, and its valency is $4590$.
Moreover, $N_G(A)/A$ acts transitively on the vertex set, and the stabilizer of a vertex is isomorphic to $N_G(S^*)/A$.
\et

\subsection{Rank $3$ property of the graph for the normalizer}
By Proposition \ref{Lgraph2}, the stabilizer of $B_\mu$ is $N_G(S^*)/A$.
In this subsection, we show that $\Gamma$ is rank $3$ for $N_G(A)/A$, that is, there are three orbits in $V\Gamma$ under the action of $N_G(S^*)/A$.

\bp\label{LG3}\begin{enumerate}
\item The subgroup $N_{G}(S^*)/A$ is transitive on the subset of $V\Gamma$ consisting of all vertices adjacent to ${{B}}_{\mu}$.
\item The subgroup $N_{G}(S^*)/A$ is transitive on the subset of $V\Gamma\setminus\{{{B}}_{\mu}\}$ consisting of all vertices not adjacent to ${{B}}_{\mu}$.
\end{enumerate}
\ep
\proof\ (1) Let $I$ be the subset of $V\Gamma$ consisting of vertices adjacent to ${{B}}_{\mu}$.
Then $\bigoplus_{{{B}}_\chi\in I}{{B}}_\chi$ is a module for $N_G(S^*)$.
By Lemmas \ref{LGr1}, $\bigoplus_{{{B}}_\chi\in I}{{B}}_\chi=(U^{\otimes 2})_2\cap B$, and by Lemma \ref{L35} it is irreducible for $N_G(S^*)$.
Hence the transitivity of $N_G(S^*)$ on $I$ follows from Lemma \ref{LG2} (1).

(2) Let $B_\chi\in V\Gamma\setminus\{B_{\mu}\}$ not adjacent to $B_{\mu}$.
Then by Lemma \ref{LGr1} (1), $B_\chi\in (M^{\otimes 2})_2$ for some $[M]\in S\setminus\{[U]\}$. 
By (\ref{ng}) and (\ref{ES:C}), $N_G(S^*)/C_G(S^*)\cong\Omega^+(10,2)$, and it is transitive on the set of non-zero singular vectors in $S=R(U)\cong \F_2^{10}$.
By Proposition \ref{POVE} (3), for $[M]\in S\setminus\{[U]\}$, $(M^{\otimes 2})_2\neq0$ if and only if $[M]$ is non-zero singular.

Let us show that $C_G(S^*)$ acts transitively on the set of all vertices in $(M^{\otimes 2})_2$ for  $M=V_{\sBW}^-$, which proves (2).
Since $\sBW$ has no roots, $M_1\cong\C\otimes_\Z\sBW$.
Hence $(M^{\otimes 2})_2=M_1^{\otimes2}\cong (\C\otimes_\Z\sBW)^{\otimes2}$ as vector spaces.
By Lemma \ref{Lfixpt} (1) and the definition of $z$, $\Ker\eta\circ\nu$ acts trivially on $(M^{\otimes 2})_2$ (see (\ref{Seq2}) and (\ref{Seq1}) for $\nu$ and $\eta$).
By $C_G(S^*)\subset C_G(z)$, $C_G(S^*)$ acts on $(M^{\otimes 2})_2$ as $\eta\circ\nu(C_G(S^*))$.
By (\ref{ES:C}), $\eta\circ\nu(\rho^{-1}(C\times C))\subset\eta\circ\nu(C_G(S^*))$.
By (\ref{Eq:C2}) and Lemma \ref{LCent1} (4), $\eta\circ\nu(\rho^{-1}(C\times C))=O_2(\Aut \sBW)*O_2(\Aut\sBW)\cong 2^{1+16}_+$.
Recall that $M_1$ is a $2^4$-dimensional irreducible module for $O_2(\Aut\sBW)\cong 2^{1+8}_+$.
By Lemma \ref{L12}, $M_1^{\otimes 2}$ is irreducible for $\eta\circ\nu(\rho^{-1}(C\times C))$.
Hence by Lemma \ref{LG2} (1), $C_G(S^*)$ acts transitively on the vertices in $(M^{\otimes 2})_2$.
\qe

Proposition \ref{LG3} also shows that the action of $N_{G}(A)/A$ on $V\Gamma$ is primitive.
Hence $N_G(S^*)/A$ is maximal in $N_{G}(A)/A$.

\bt\label{TM1} The graph $\Gamma$ is rank $3$ for $N_{G}(A)/A$.
Moreover, $N_{G}(A)$ is generated by $N_{G}(S^*)$ and $C_{G}(z)$.
\et

\subsection{Shape of the normalizer}
In this subsection, we determine the shape of $N_G(A)$ by using the classification of permutation groups of rank $3$.

As a corollary of Theorem \ref{TM1}, we have the following.

\bc The group $N_G(A)/A$ is a primitive permutation group of rank $3$ of degree $139503=3\cdot7^2\cdot13\cdot73$.
\ec

Now, we recall that any primitive permutation group $Q$ of rank $3$ of finite degree $n$ satisfies one of the following (see \cite{Li} and references therein):
\begin{enumerate}[(i)]
\item $Q$ is a subgroup of $E_0\wr \Z_2$ and $Q$ contains $E\times E$ as a normal subgroup, where $E_0$ is a $2$-transitive group of degree $n_0$ and the socle $E$ of $E_0$ is a simple group, and $n=n_0^2$;
\item The socle of $Q$, the product of all the minimal normal subgroups of $Q$, is (non-abelian) simple;
\item $Q$ is an affine group, that is, the socle of $Q$ is a vector space $W$, where $W\cong \Z_p^d$ for some prime $p$ and $n=p^d$.
\end{enumerate}

By the degree, $N_G(A)/A$ satisfies (ii).
Since the stabilizer in $N_G(A)/A$ of a vertex in $\Gamma$ has a shape $2^{16}:\Omega^+(10,2)$ and its index is $139503$, the order of $N_G(A)/A$ is $2^{36}\cdot3^6\cdot5^2\cdot7^3\cdot13\cdot17\cdot 31\cdot 73$; the socle of $N_G(A)/A$ is isomorphic to neither an alternating group nor a sporadic group (cf.\ \cite{ATLAS}); the primitive permutation groups of rank $3$ under (ii) were determined in \cite{KL} when the socle is a classical group, in \cite{LS} when the socle is a non-classical simple group of Lie type; $N_G(A)/A$ must be isomorphic to $E_6(2)$.

\bp\label{PE6} The normalizer $N_G(A)$ has the shape $2^{27}.E_6(2)$.
\ep

By Theorem \ref{TM1} and Proposition \ref{PE6} (cf.\ \cite{LS}), we obtain the following corollary:

\bc\label{MC} The graph $\Gamma$ is isomorphic to the rank $3$ graph on $E_6(2)/(2^{16}: \Omega^+(10,2))$ of valency $4590$, and its automorphism group is isomorphic to $E_6(2)$.
\ec

\subsection{Full automorphism group}

In this subsection, we determine the full automorphism group $G$ of the VOA $V$.

First, we show the finiteness of $G$.
The following is a modification of \cite[Theorem (iv)]{Ti}.

\bp\label{P5-1-1} Let $\mathcal{H}$ be an algebraic group.
Assume that $\mathcal{H}$ contains commuting involutions $z_1,z_2$ such that $C_\mathcal{H}(z_1)$, $C_\mathcal{H}(z_2)$ and $C_\mathcal{H}(z_1z_2)$ are finite.
Then $\mathcal{H}$ is finite.
\ep
\proof\ Let $\mathcal{H}^0$ be the connected component of the identity in $\mathcal{H}$.
Clearly, $\mathcal{H}^0$ is normal in $\mathcal{H}$.
Let $a$ be an involution in $\mathcal{H}$ such that $C_\mathcal{H}(a)$ is finite.
Then the involution $g\mapsto a^{-1}ga$ in $\Aut \mathcal{H}^0$ has only finitely many fixed points.
Let $L(\mathcal{H}^0)$ be the Lie algebra of $\mathcal{H}^0$.
Then the canonical group homomorphism $\phi$ from $\Aut(\mathcal{H}^0)$ to $\Aut(L(\mathcal{H}^0))$ is injective.

We now suppose that $\mathcal{H}^0\neq\{1\}$.
Since $C_{\mathcal{H}}(a)$ is finite, $\phi(a)$ is an involution in $\Aut L(\mathcal{H}^0)$ whose only fixed point is $0$.
This implies that $\phi(a)=-1$.
In particular, $\phi(z_1)=\phi(z_2)=\phi(z_1z_2)=-1$, which contradicts $\phi(z_1)\phi(z_2)=\phi(z_1z_2)$.
Hence $\mathcal{H}^0=\{1\}$, and $\mathcal{H}$ is discrete.
Therefore $\mathcal{H}$ is finite since it is an algebraic group.
\qe

By Proposition \ref{PCen} $C_G(z)$ is finite.
We now check that $G$ satisfies the hypothesis of the proposition above.

\bl\label{L5-1-1} The abelian subgroup $S^*$ of $G$ contains commuting involutions $z_1,z_2$ such that $z_1$, $z_2$ and $z_1z_2$ are conjugate to $z$.
In particular, $C_G(z_1)$, $C_G(z_2)$ and $C_G(z_1z_2)$ are finite.
\el
\proof\ By the natural correspondence between $S^*$ and $S$ with respect to the symplectic form, $z$ corresponds to the non-zero singular vector $[V_{\sBW}^-]\in S$.
The group $N_G(S^*)$ acts on $S^*$ as a natural module for $\Omega^+(10,2)$ by conjugation, and it is transitive on the set of non-zero singular vectors in $S$.
Hence all involutions in $S^*$ corresponding to non-zero singular vectors in $S$ are conjugate to $z$.
Since a $10$-dimensional non-singular quadratic space over $\F_2$ contains a $2$-dimensional totally singular subspace, we obtain this lemma.
\qe

It was shown in \cite[Theorem 2.4]{DG2} that the automorphism group of a finitely generated VOA is an algebraic group.
By Lemma \ref{Lwt2}, $V$ is finitely generated, and hence $G$ is an algebraic group.
Thus Proposition \ref{P5-1-1} and Lemma \ref{L5-1-1} imply the following:

\bp The automorphism group $G$ of $V$ is finite.
\ep

Next, we will show that $A$ is normal in $G$.

\bl\label{LS1} Let $P$ be a Sylow $2$-subgroup of $C_G(z)$.
Then $Z(P)=\langle z\rangle$, and $P$ is also a Sylow $2$-subgroup of $G$.
\el
\proof\ It suffices to consider a specified Sylow $2$-subgroup $P\cong 2^{1+32}_+.(Q_1.Q_2)$ of $C_G(z)$, where $Q_1=O_2(\Aut\BW)$ and $Q_2$ is a Sylow $2$-subgroup of $\Aut\BW/O_2(\Aut\BW)$.
Then $Q_2$ and $Q_1$ act faithfully on $Q_1$ and $2^{1+32}/\langle z\rangle\cong \BW/2\BW$, respectively.
Hence $Z(P)\subset Z(2^{1+32}_+)=\langle z\rangle$, and $Z(P)=\langle z\rangle$. 
Let $P_0$ be a Sylow $2$-subgroup of $G$ containing $P$.
Then $1\neq Z(P_0)\subset Z(P)=\langle z\rangle$, and hence $Z(P_0)=Z(P)$.
Thus $P_0\subset C_G(z)$, and $P_0=P$.
\qed

\bp\label{PNormal} Let $F$ be a minimal normal subgroup of $G$.
Then $F$ is abelian, and $F=A$.
\ep
\proof\ Set $J=F\cap N_G(A)$.
By Proposition \ref{PAirr} $A$ is irreducible as an $N_G(A)$-module, and by Proposition \ref{PE6} $N_G(A)/A$ is simple.
Hence $J$ must be $\{1\}$, $A$ or $N_G(A)$.

First, we prove $A\subset F$.
Suppose $A\not\subset F$.
Then $J=\{1\}$ and $C_G(z)\cap F=\{1\}$.
Let $x,y\in F$ with $x^zx^{-1}=y^zy^{-1}$.
Then $(y^{-1}x)^z=(y^{-1})^zx^z=y^{-1}x$.
By $C_G(z)\cap F=\{1\}$, we have $y^{-1}x=1$, and $y=x$.
Moreover, if $x\neq 1$ then $x^zx^{-1}\neq 1$.
Hence $F=\{x^zx^{-1}\mid x\in F\}$, and $f^z=f^{-1}$ for all $f\in F$.
By Lemma \ref{L5-1-1}, there exist involutions $z_1,z_2\in S^*$ such that $z_1$, $z_2$ and $z_1z_2$ are conjugate to $z$.
The argument above shows that $z_1$, $z_2$, $z_1z_2$ induce the same automorphism of $F$ of order $2$, which is a contradiction.

Next, we suppose that $F$ is not abelian.
Then $F$ is a direct product of isomorphic (non-abelian) simple groups, that is, $F=\Pi_{i=1}^lF_i$, $F_i\cong F_j$ and $F_i$ is simple for all $i,j$; $A=\Pi_{i=1}^l(A\cap F_i)$ by $N_G(A)\cap F=A$.
By $C_N(A)=A$, $A\cap F_i\neq \{1\}$ for all $i$.
Hence $N_G(A)/A$ acts on $\{A\cap F_i\}$ as a permutation group of degree $d\le 27$.
Since $N_G(A)/A$ is simple and $73\mid |E_6(2)|$, $d$ must be $1$; $F$ is simple.

We consider the following cases: (1) $J=A$; (2) $J=N_G(A)$.

(1) Suppose $J=A$.
By Lemma \ref{LS1}, there exists a Sylow $2$-subgroup $P$ of $G$ such that $A\subset P\subset C_G(z)$.
Then $P\cap F=A$ is a Sylow $2$-subgroup of $F$.
By $J=A$, we have $N_F(A)=C_F(A)=A$.
By \cite[Theorem 4.3 in Chapter 7]{Go}, $F$ has a normal complement of $A$, which contradicts the simplicity of $F$.

(2) Suppose $J=N_G(A)$.
This implies $N_G(A)\subset F$.
It follows from $C_G(z)\subset N_G(A)$ that $C_F(z)=C_G(z)$.
Let $P$ be a Sylow $2$-subgroup of $C_G(z)$.
Then by Lemma \ref{LS1} $P$ is also a Sylow $2$-subgroup of $F$.
By Proposition \ref{PCen}, $|P|=2^{63}$. 
By $N_G(A)\subset F$ and Proposition \ref{PE6}, $|F|$ is divisible by $73$.
We now refer the classification of finite simple groups.
If $F$ is an alternating group, then by $|P|=2^{63}$ it isomorphic to $A_{66}$ or $A_{67}$, which contradicts that $|F|$ is divisible by $73$.
By \cite{ATLAS} and $|P|=2^{63}$, $F$ is not a sporadic simple group.
According to \cite{AS} and \cite[Section 4.5]{GLS}, $C_G(z)\cong 2^{1+32}.2^{10}.\Omega^+(10,2)$ can not be the centralizer of an involution in any finite simple group of Lie type.
Thus (2) is impossible.

Therefore $F$ is abelian and $A\subset F$.
By Proposition \ref{CG2}, $A\subset F\subset C_G(A)=A$, and hence $F=A$.
\qe

Combining Propositions \ref{PE6} and \ref{PNormal}, we obtain the following main theorem:

\bt\label{MT} The automorphism group of the $\Z_2$-orbifold VOA $\tilde{V}_{\BW}$ has the shape $2^{27}.E_6(2)$.
\et

\paragraph{Acknowledgments.} 
The author thanks Professor Robert L. Griess Jr. for helpful comments, Professor Masahiko Miyamoto for sending the unpublished manuscript \cite{Mi}, and Professor Ching Hung Lam for carefully reading an early version of this article.
He also thanks the referee for many valuable suggestions.
\begin{small}

\end{small}


\begin{thebibliography}{100000}

\bibitem[AD04]{AD}
T.\ Abe and C.\ Dong, Classification of irreducible modules for the vertex operator algebra $V\sp +\sb L$: general case. {\it J. Algebra} {\bf 273} (2004), 657--685

\bibitem[ADL05]{ADL}
T.\ Abe, C.\ Dong and H.\ Li, Fusion rules for the vertex operator algebras $M(1)^+$ and $V_L^+$, {\it Comm.\ Math. Phys.} {\bf 253} (2005), 171--219.

\bibitem[AS76]{AS}
M.\ Aschbacher and G.M.\ Seitz, Involutions in Chevalley groups over fields of even order, {\it Nagoya Math. J.} {\bf 63} (1976), 1--91. 

\bibitem[BW59]{BW}
E.S.\ Barnes and G.E.\ Wall, Some extreme forms defined in terms of Abelian groups, {\it J.\ Austral.\ Math.\ Soc.} {\bf 1} (1959), 47--63

\bibitem[Bo86]{Bo}
R.E.\ Borcherds, Vertex algebras, Kac-Moody algebras, and the Monster, {\it Proc.\ Nat'l.\ Acad.\ Sci.\ U.S.A.} {\bf 83} (1986), 3068--3071.

\bibitem[Bo92]{Bo2}
R.E. Borcherds, Monstrous moonshine and monstrous Lie superalgebras, {\it Invent. Math.} {\bf 109} (1992), 405--444.

\bibitem[CCNPW85]{ATLAS}
J.H.\ Conway, R.T.\ Curtis, S.P.\ Norton, R.A.\ Parker and R.A.\ Wilson, Atlas of finite groups, Oxford University Press, Oxford, 1985.

\bibitem[CS99]{CS}
J.H.\ Conway and N.J.A. Sloane, Sphere packings, lattices and groups, 3rd Edition, Springer, New York, 1999.


\bibitem[DGM96]{DGM}
L.\ Dolan, P.\ Goddard and P.\ Montague, Conformal field theories, representations and lattice constructions, {\it Comm. Math. Phys.} {\bf 179} (1996), 61--120.

\bibitem[Do93]{Do}
C.\ Dong, Vertex algebras associated with even lattices, {\it J. Algebra} {\bf 161} (1993), 245--265.

\bibitem[DG98]{DG0}
C.\ Dong and R.L.\ Griess, Rank one lattice type vertex operator algebras and their automorphism groups, {\it J.\ Algebra} {\bf 208} (1998), 262--275.

\bibitem[DG02]{DG2}
C.\ Dong and R.L.\ Griess, Automorphism groups and derivation algebras of finitely generated vertex operator algebras.  {\it Michigan Math.\ J.} {\bf 50} (2002), 227--239.

\bibitem[DG05]{DG}
C.\ Dong and R.L.\ Griess, The rank two lattice type vertex operator algebras $V_L^+$ and their automorphism groups,  {\it Michigan Math.\ J.} {\bf 53} (2005), 691--715.

\bibitem[DGH98]{DGH}
C.\ Dong, R.L.\ Griess, and G.\ H$\ddot{\rm o}$hn, Framed vertex operator algebras, codes and Moonshine module, {\it Comm.\ Math.\ Phys.} {\bf 193} (1998), 407--448.

\bibitem[DLM00]{DLM}
C.\ Dong, H.\ Li, and G.\ Mason, Modular-invariance of trace functions in orbifold theory and generalized Moonshine, {\it Comm. Math. Phys.} {\bf 214} (2000), 1--56.

\bibitem[DM04]{DM}
C.\ Dong and G.\ Mason, Rational vertex operator algebras and the effective central charge, {\it Int. Math. Res. Not.} (2004), 2989--3008.

\bibitem[DN99]{DN}
C.\ Dong and K.\ Nagatomo, Automorphism groups and Twisted modules for lattice Vertex operator algebras, {\it Comtemp.\ Math.} {\bf 248} (1999), 117--133

\bibitem[FHL93]{FHL}
I.\ Frenkel, Y.\ Huang and J.\ Lepowsky, On axiomatic approaches to vertex operator algebras and modules,  Mem. Amer. Math. Soc. {\bf 104} 1993.

\bibitem[FLM88]{FLM}
I.\ Frenkel, J.\ Lepowsky and A.\ Meurman, Vertex operator algebras and the Monster, Pure and Appl.\ Math., Vol.134, Academic Press, Boston, 1988.

\bibitem[Go80]{Go}
D.\ Gorenstein, Finite groups, Chelsea Publishing Co., New York, 1980.

\bibitem[GLS98]{GLS}
D.\ Gorenstein, R.\ Lyons and R.\ Solomon, The classification of the finite simple groups, Number 3. Mathematical Surveys and Monographs, {\bf 40.3.}, American Mathematical Society, Providence, 1998.

\bibitem[Gr82]{Gr}
R.L.\ Griess, The friendly giant. {\it Invent. Math.} {\bf 69} (1982) 1--102. 


\bibitem[Gr98]{Gr1}
R.L.\ Griess, A vertex operator algebra related to $E\sb 8$ with automorphism group ${\rm O}\sp +(10,2)$, {\it Ohio State Univ. Math. Res. Inst. Publ.} {\bf 7} (1998), 43--58.

\bibitem[Gr05]{Gr3}
R.L.\ Griess,  Pieces of $2\sp d$: existence and uniqueness for Barnes-Wall and Ypsilanti lattices, {\it Adv. Math.} {\bf 196} (2005), 147--192, Corrections and additions, {\bf 211} (2007), 819--824.

\bibitem[H\"o95]{HoPhd}
G. H\"ohn, Selbstduale Vertexoperatorsuperalgebren und das Babymonster, PhD thesis, Universit\"at Bonn, 1995, Bonner Math. Schriften {\bf 286}, 1996.

\bibitem[H\"{o}08]{Ho}
G.\ H\"{o}hn, Conformal designs based on vertex operator algebras, {\it Adv. Math.} {\bf 217} (2008), 2301--2335.

\bibitem[HL95]{HL}
Y.-Z\ Huang and J.\ Lepowsky, A theory of tensor products for module categories for a vertex operator algebra. III, {\it J. Pure Appl. Algebra} {\bf 100} (1995), 141--171.

\bibitem[KL82]{KL}
W.M.\ Kantor and R.A.\ Liebler, The rank $3$ permutation representations of the finite classical groups, {\it Trans. Amer. Math. Soc.} {\bf 271}  (1982), 1--71.

\bibitem[LSY07]{LSY}
C.\ Lam, S.\ Sakuma and H.\ Yamauchi, Ising vectors and automorphism groups of commutant subalgebras related to root systems, {\it Math.\ Z.} {\bf 255} (2007), 597--626. 

\bibitem[Li87]{Li}
M.W.\ Liebeck, The affine permutation groups of rank three, {\it Proc. London Math. Soc.} (3) {\bf 54} (1987), 477--516.

\bibitem[LS86]{LS}
M.W.\ Liebeck and J.\ Saxl, The finite primitive permutation groups of rank three, {\it Bull. London Math. Soc.} {\bf 18} (1986), 165--172.

\bibitem[Ma01]{Ma}
A.\ Matsuo, Norton's trace formulae for the Griess algebra of a vertex operator algebra with larger symmetry, {\it Comm. Math. Phys.} {\bf 224} (2001), 565--591.

\bibitem[MM00]{MM}
A.\ Matsuo and M.\ Matsuo, The automorphism group of the Hamming code vertex operator algebra, {\it J.\ Algebra} {\bf 228} (2000), 204--226.

\bibitem[Mi96]{Mi0} M.\ Miyamoto, Griess algebras and conformal vectors in vertex
operator algebras, {\it J. Algebra} {\bf 179} (1996), 523--548.

\bibitem[Mi]{Mi}
M.\ Miyamoto, Automorphism groups of $\Z_2$-orbifold VOAs, unpublished paper, 1996.

\bibitem[Mi04]{Mi3}
M. Miyamoto, A new construction of the Moonshine vertex operator algebra over the real number field, {\it Ann.\ of Math.} {\bf 159} (2004), 535--596.

\bibitem[Sh04]{Sh2}
H.\ Shimakura, The automorphism group of the vertex operator algebra $V_L^+$ for an even lattice $L$ without roots, {\it J.\ Algebra} {\bf 280} (2004), 29--57.

\bibitem[Sh06]{Sh3}
H.\ Shimakura, The automorphism groups of the vertex operator algebras $V\sp +\sb L$: general case, {\it Math.\ Z.} {\bf 252} (2006), 849--862.

\bibitem[Sh07]{Sh4}
H.\ Shimakura, Lifts of automorphisms of vertex operator algebras in simple current extensions, {\it Math.\ Z.} {\bf 256} (2007), 491--508.

\bibitem[Sh11]{Sh6}
H.\ Shimakura, An $E_8$-approach to the moonshine vertex operator algebra, {\it J. London Math. Soc.} {\bf 83} (2011), 493--516.

\bibitem[Sh12]{Sh10}
H.\ Shimakura, Classification of Ising vectors in the vertex operator algebra $V_L^+$, {\it Pacific J. Math.} {\bf 258} (2012), 487--495. 


\bibitem[Ti84]{Ti}
J.\ Tits, On R. Griess' "friendly giant", {\it Invent.\ Math.} {\bf 78} (1984), 491--499.

\end{thebibliography}
\end{document}